\documentclass{amsart}
\theoremstyle{plain}
\newtheorem{theorem}{Theorem}

\newtheorem{proposition}[theorem]{Proposition}
\theoremstyle{definition}

\theoremstyle{remark}
\newtheorem{remark}[theorem]{Remark}

\title{On Vandermonde determinants via  $n$-determinants}
\author{Milan Janji\'c}
\date{\today}
\address{Department of Mathematics and Informatics\\
 University of Banja Luka \\
Republic of Srpska, BA}
\begin{document}
\maketitle
\begin{center}Department for Mathematics and Informatics, University of Banja
Luka\\Republic of Srpska, Bosnia and Herzegovina\end{center}

\begin{abstract}
We use earlier defined notion of $n$- determinant to investigate sub-determinants of an extended Vandermonde matrix. Firstly, we demonstrate our method on a number of particular cases. Then we prove that all these results may be stated in terms of Schur's polynomials. In our main result, we prove that Schur polynomials are equal to minors of a fixed matrix, which entries are formed of elementary symmetric polynomials. Such a formula is known as the  second Jaccobi-Trudi identity.
\end{abstract}

\section{Introduction}
Let $X=\{x_1,x_2,\ldots,x_n\}$ be a set of variables.
 We consider the following extended Vandermonde matrix:
 \begin{displaymath}
V_{n,n+r}=
\begin{pmatrix}1&x_1&x_1^2&\cdots&x_1^{n}&\cdots&x_1^{n+r-1}\\
1&x_2&x_2^2&\cdots&x_2^{n}&\cdots&x_2^{n+r-1}\\
\vdots&\vdots&\vdots&\cdots&\vdots&\cdots&\vdots\\
1&x_n&x_n^2&\cdots&x_n^n&\cdots&x_n^{n+r-1}
\end{pmatrix}.\end{displaymath}
First $n$ columns of $V_{n,n+r}$ form the standard Vandermonde matrix
\begin{displaymath}V_n=\begin{pmatrix}1&x_1&\cdots&x_1^{n-1}\\
1&x_2&\cdots&x_2^{n-1}\\
\vdots&\vdots&\cdots&\vdots\\
1&x_n&\cdots&x_n^{n-1}
\end{pmatrix}.
\end{displaymath}
We  denote by $K_i,(i=1,2,\ldots,n+r)$  columns of $V_{n,n+r}$.
 By
$\overline{e}_i(X),(i=0,1,\ldots,n)$ are denoted elementary symmetric polynomial of
$x_1,x_2,\ldots,x_n$. Hence, $\overline{e}_0(X)=1,\overline{e}_1(X)=
x_1+x_2+\cdots+x_n,\ldots,\overline{e}_n(X)=x_1\cdot x_2\cdots x_n$.

We next consider the following polynomial
$p_n(x)=\prod_{i=1}^n(x-x_i)$. Expanding the right-hand side we obtain
\begin{displaymath}p_n(x)=\sum_{i=0}^n(-1)^{n-i}\cdot \overline{e}_{n-i}(X)\cdot
x^{i}\end{displaymath}
From $p_n(x_k)=0,(k=1,\ldots,n)$, we obtain

\begin{displaymath}x_k^n=\sum_{i=1}^n(-1)^{n-i}\overline{e}_{n-i+1}(X)\cdot
x_k^{i-1}.
\end{displaymath}
Denoting $(-1)^{n-i}\overline{e}_{n-i+1}(X)=e_{n-i+1}(X),(i=1,\ldots,n)$   implies
\begin{displaymath}
K_{n+1}=\sum_{i=1}^ne_{n-i+1}(X)\cdot K_{i}.
\end{displaymath}.
In the same way, for each $p=1,2,\ldots,r$, we obtain
\begin{equation}\label{kki}
K_{n+p}=\sum_{i=1}^ne_{n-i+1}(X)\cdot K_{i+p-1}.
\end{equation}.
\begin{remark} We may consider that  each column $K_{n+k}$ is, in fact,
 a linear combinations of all preceding columns of $V_{n,n+r}$,
  by taking coefficient $0$ for columns that do not appear in the equation
  (\ref{kki}).
\end{remark}

We next consider the following matrix of order $(r+n-1)\times r$.
\begin{displaymath}P_{n+r-1,r}=\begin{pmatrix}
e_n&0&0&\cdots&0&0&0\\
e_{n-1}&e_n&0&\cdots&0&0&0\\
\vdots&\vdots&\vdots&\cdots&\vdots&\vdots&\vdots\\
e_1&e_2&e_3&\cdots&0&0&0\\
-1&e_1&e_2&\cdots&0&0&0\\
\vdots&\vdots&\vdots&\cdots&\vdots&\vdots&\vdots\\
0&0&0&\cdots&-1&e_1&e_2\\
0&0&0&\cdots&0&-1&e_1
\end{pmatrix}.\end{displaymath}

Assume that
\begin{displaymath}\{1,2,\ldots,n+r-1\}=\{i_1,i_2,\ldots,i_{n-1}\}
\cup\{j_1,j_2,\ldots,j_{r}\}.\end{displaymath}
Clearly, the union on the right-hand side must be disjoint.
We  denote by
\begin{displaymath}M=M(i_1,i_2,\ldots,i_{n-1},n+r)\end{displaymath} the sub-matrix of
$V_{n,n+r}$ lying in columns
$i_1,i_2,\ldots,i_{n-1},n+r$ of $V_{n,n+r}$. We  define
 \begin{displaymath}{\rm sgn
 }M=(-1)^{\frac{n(n-1)}{2}+i_1+\cdots+i_{n-1}}.\end{displaymath}
 It is easy to see that the equation
 ${\rm sgn }M=(-1)^{nr+\frac{r(r-1)}{2}+j_1+\cdots+j_{r}}$ also holds.

 We next denote by $Q_r=Q_r({j_1,j_2,\ldots,j_r})$ the sub-matrix of order $r$ of
 $P$ lying in rows $j_1,j_2,\ldots,j_r$ of $P$.

 In our paper \cite{jan}, the following formula is proved:
 \begin{equation}\label{fun}
 \frac{\begin{vmatrix}x_1^{i_1-1}&\cdots&x_1^{i_{n-1}-1}&x_1^{n+r-1}\\
x_2^{i_1-1}&\cdots&x_2^{i_{n-1}-1}&x_2^{n+r-1}\\
\vdots&\cdots&\vdots&\vdots\\
x_n^{i_1-1}&\cdots&x_n^{i_{n-1}-1}&x_n^{n+r-1}
\end{vmatrix}}{\begin{vmatrix}1&x_1&\cdots&x_1^{n-1}\\
1&x_2&\cdots&x_2^{n-1}\\
\vdots&\vdots&\cdots&\vdots\\
1&x_n&\cdots&x_n^{n-1}
\end{vmatrix}}={\rm sgn }M\cdot\det Q_r(j_1,\ldots,j_r).
 \end{equation}

We start with
examples of this result.
\section{Two introducing examples}  In this part, we describe two extreme case,
when $Q_r(j_1,\ldots,j_r)$ lies either in the first $r$ rows of $P$ or in
 the last $r$.
 We firstly take $i_1=r+1,i_2=r+2,\ldots,i_{n-1}=r+n-1$.
  Then $Q_r$ lies in the first $r$ rows of $P$. Hence,
  it is a lower triangular matrix having all
  diagonal elements equal to $e_n$. Also, ${\rm sgn }M=(-1)^{(n-1)\cdot r}$, so that
  \begin{proposition}The following equation holds
 \begin{displaymath}\frac{
\begin{vmatrix}x_1^{r}&x_1^{r+1}&\cdots&x_1^{r+n-2}&x_1^{r+n-1}\\
x_2^{r}&x_2^{r+1}&\cdots&x_2^{r+n-2}&x_2^{r+n-1}\\
\vdots&\vdots&\cdots&\vdots&\vdots\\
x_n^{r}&x_n^{r+1}&\cdots&x_n^{r+n-2}&x_n^{r+n-1}
\end{vmatrix}}{\begin{vmatrix}1&x_1&\cdots&x_1^{n-2}&x_1^{n-1}\\
1&x_2&\cdots&x_2^{n-2}&x_2^{n-1}\\
\vdots&\vdots&\cdots&\vdots&\vdots\\
1&x_n&\cdots&x_n^{n-2}&x_n^{n-1}
\end{vmatrix}}=(-1)^{(n-1)\cdot r}\cdot e_n^r(X).
\end{displaymath}
Extracting $x_1^rx_2^r\cdots x_n^r$ from nominator of the fraction, the equation becomes obvious.
\end{proposition}
\begin{proposition} Assume that $i_t=t,(t=1,2,\ldots,n-1)$. We obviously have
${\rm sgn } M=1$.  For $r>n$, we have
\begin{displaymath}
\frac{\begin{vmatrix}1&x_1&\cdots&x_1^{n-2}&x_1^{n+r-1}\\
1&x_2&\cdots&x_2^{n-2}&x_2^{n+r-1}\\
\vdots&\vdots&\cdots&\vdots&\vdots\\
1&x_n&\cdots&x_n^{n-2}&x_n^{n+r-1}
\end{vmatrix}}{\begin{vmatrix}1&x_1&\cdots&x_1^{n-2}&x_1^{n-1}\\
1&x_2&\cdots&x_2^{n-2}&x_2^{n-1}\\
\vdots&\vdots&\cdots&\vdots&\vdots\\
1&x_n&\cdots&x_n^{n-2}&x_n^{n-1}
\end{vmatrix}}=\begin{vmatrix}
e_1&e_2&\cdots&0&0\\
-1&e_1&\cdots&0&0\\
\vdots&\vdots&\cdots&\vdots&\vdots\\
0&0&\cdots&e_1&e_2\\
0&0&\cdots&-1&e_1
\end{vmatrix}.
\end{displaymath}
\end{proposition}

\section{Case $r=1$}
In the case $r=1$, the matrix $P$ has the following form
\begin{displaymath}P=\begin{pmatrix}
e_n(X)\\
e_{n-1}(X)\\
\vdots\\
e_1(X)\\
-1
\end{pmatrix},\end{displaymath}
  If $A_1,A_2,\ldots,A_n,A_{n+1}$ are columns of $V_{n,n+1}$, then we have
\begin{displaymath}A_{n+1}=\sum_{i=1}^ne_{n-i+1}(X)\cdot A_{n-i+1}.
\end{displaymath}
We next have  ${\rm sgn} M=\frac{n^2-n+2}{2}$. Hence,
\begin{proposition}
We have
\begin{equation}\label{eq1}
\frac{
\begin{vmatrix}
1&x_1&\cdots&x_{1}^{j-2}&x_{1}^{j}&\cdots&x_{1}^{n}\\
1&x_2&\cdots&x_{1}^{j-2}&x_{2}^{j}&\cdots&x_{2}^{n}\\
\vdots&\vdots&\cdots&\vdots&\vdots&\cdots&\vdots\\
1&x_n&\cdots&x_{n}^{j-2}&x_{n}^{j}&\cdots&x_{n}^{n}
\end{vmatrix}}{\begin{vmatrix}1&x_1&\cdots&x_1^{n-2}&x_1^{n-1}\\
1&x_2&\cdots&x_2^{n-2}&x_2^{n-1}\\
\vdots&\vdots&\cdots&\vdots&\vdots\\
1&x_n&\cdots&x_n^{n-2}&x_n^{n-1}
\end{vmatrix}}=(-1)^{\frac{n^2-n+2}{2}}\cdot
 e_{n-j+1}(X).\end{equation}
\end{proposition}

\section{Case $r=2$}
In this case, we have
\begin{displaymath}
V_{n,n+2}=\begin{pmatrix}1&x_1&x_1^2&\cdots&x_1^{n}&x_1^{n+1}\\
1&x_2&x_2^2&\cdots&x_2^{n}&x_2^{n+1}\\
\vdots&\vdots&\vdots&\cdots&\vdots&\vdots\\
1&x_{n-1}&x_{n-1}^2&\cdots&x_{n-1}^{n}&x_{n-1}^{n+1}\\
1&x_n&x_n^2&\cdots&x_n^{n}&x_n^{n+1}.
\end{pmatrix}
\end{displaymath}
We denote by $\hat{x}_i$ deleted elements of $V_{n,n+2}$. If $i$ and $j$ are indices
of deleted columns, we have
\begin{displaymath}M=
\begin{pmatrix}1&x_1&\cdots&\hat{x}_{1}^{i-1}&
\cdots&\hat{x}_{1}^{j-1}&\cdots&x_1^{n+1}\\
1&x_2&\cdots&\hat{x}_{2}^{i-1}&
\cdots&\hat{x}_{2}^{j-1}&\cdots&x_2^{n+1}\\
\vdots&\vdots&\cdots&\vdots&\vdots\\
1&x_n&\cdots&\hat{x}_{n}^{i-1}&
\cdots&\hat{x}_{n}^{j-1}&\cdots&x_n^{n+1}
\end{pmatrix}.
\end{displaymath}
In this case the matrix $P$ is of order $(n+1)\times 2$, and $Q$ is
a $2\times 2$ sub-matrix of $P$. We consider the case when
 $Q_2=\begin{pmatrix}e_i(X)&e_{i+1}(X)
\\e_j(X)&e_{j+1}(X)\end{pmatrix}$.
Next, we have ${\rm sgn} M=(-1)^{nr+\frac{r(r-1)}{2}+i+j}$. We thus obtain
\begin{proposition} The following formula holds
\begin{displaymath}
\frac{\begin{vmatrix}1&\cdots&\hat{x}_{1}^{i-1}&
\cdots&\hat{x}_{1}^{j-1}&\cdots&x_1^{n+1}\\
1&\cdots&\hat{x}_{2}^{i-1}&
\cdots&\hat{x}_{2}^{j-1}&\cdots&x_2^{n+1}\\
\vdots&\cdots&\vdots&\vdots\\
1&\cdots&\hat{x}_{n}^{i-1}&
\cdots&\hat{x}_{n}^{j-1}&\cdots&x_n^{n+1}
\end{vmatrix}}{\begin{vmatrix}1&x_1&\cdots&x_1^{n-2}&x_1^{n-1}\\
1&x_2&\cdots&x_2^{n-2}&x_2^{n-1}\\
\vdots&\vdots&\cdots&\vdots&\vdots\\
1&x_n&\cdots&x_n^{n-2}&x_n^{n-1}
\end{vmatrix}
}=(-1)^{nr+\frac{r(r-1)}{2}+i+j}\cdot \begin{vmatrix}e_i(X)&e_j(X)\\
 e_{i+1}(X)&e_{j+1}(X)\end{vmatrix}
.
\end{displaymath}
\end{proposition}

\section{Case $i_1=1,i_{2}=2,\ldots,i_{n-1}=n-1$}
In this case, we have
\begin{displaymath}
M=\begin{pmatrix}1&x_1&\cdots&x_1^{n-2}&x_1^{n+r}\\
1&x_2&\cdots&x_2^{n-2}&x_2^{n+r}\\
\vdots&\vdots&\cdots&\vdots&\vdots\\
1&x_n&\cdots&x_n^{r+n-2}&x_n^{n+r}
\end{pmatrix},
\end{displaymath}
and ${\rm sgn }M=1$. Next, we have
\begin{displaymath}Q_r=\begin{pmatrix}
e_1(X)&e_2(X)&\cdots&e_n(X)&0&\cdots&0&0&0\\
-1&\overline{e}_1(X)&\cdots&e_{n-1}(X)&e_n(X)&\cdots&0&0&0\\
\vdots&\vdots&\cdots&\vdots&\vdots&\cdots&\vdots&\vdots&\vdots\\
0&0&\cdots&0&0&\cdots&-1&e_1(X)&e_2(X)\\
0&0&\cdots&0&0&\cdots&0&-1&\overline{e}_1(X)
\end{pmatrix}.\end{displaymath}
We see that $Q_r$ is a Hessenberg matrix, and also a Toeplitz matrix.
It is easy to see that $\det Q_r$ satisfies  the following recurrence:
\begin{displaymath}
\det Q_r=e_1(X)\cdot Q_{r-1}+e_2(X)\cdot \det Q_{r-2}(X)+\cdots+
e_n(X)\cdot \det Q_{r-n},(r>n).
\end{displaymath}
Being obviously ${\rm sgn }M=1$,
we obtain
\begin{proposition}The following formula holds
\begin{displaymath}
\frac{
\begin{vmatrix}1&x_1&\cdots&x_1^{n-2}&x_1^{n+r}\\
1&x_2&\cdots&x_2^{n-2}&x_2^{n+r}\\
\vdots&\vdots&\cdots&\vdots&\vdots\\
1&x_n&\cdots&x_n^{r+n-2}&x_n^{n+r}
\end{vmatrix}}{\begin{vmatrix}1&x_1&\cdots&x_1^{n-2}&x_1^{n-1}\\
1&x_2&\cdots&x_2^{n-2}&x_2^{n-1}\\
\vdots&\vdots&\cdots&\vdots&\vdots\\
1&x_n&\cdots&x_n^{n-2}&x_n^{n-1}
\end{vmatrix}}=
\det Q_r.
\end{displaymath}
\end{proposition}
\section{A relation with Schur polynomials}
In this section, we relate  obtained results with the Schur polynomials.
We  transform $M$ in the following way: Firstly,  we interchange rows
$(1,n),(2,n-1),\ldots$. For this we need $\lfloor\frac{n}{2}\rfloor$ transposition. Hence
 \begin{displaymath}
N=\begin{pmatrix}x_1^{n+r-1}&x_2^{n+r-1}&\cdots&x_n^{n+r-1}\\
x_1^{i_{n-1}-1}&x_2^{i_{n-1}-1}&\cdots&x_n^{i_{n-1}-1}\\
\vdots&\vdots&\cdots&\vdots\\
x_1^{i_1-1}&x_2^{i_1-1}&\cdots&x_n^{i_{1}-1}
\end{pmatrix},
\end{displaymath}
and  $\det N=(-1)^{\lfloor\frac{n}{2}\rfloor}\cdot \det M$.
By denoting $\lambda_1=r,\lambda_2=i_{n-1}-n+1,\ldots,\lambda_n=i_1-1$, we obtain
\begin{displaymath}
 N=\begin{pmatrix}
x_1^{\lambda_1+n-1}&x_2^{\lambda_1+n-1}&\cdots&x_n^{\lambda_1+n-1}\\
x_1^{\lambda_2+n-2}&x_2^{\lambda_2+n-2}&\cdots&x_n^{\lambda_2+n-2}\\
\vdots&\vdots&\cdots&\vdots\\
x_1^{\lambda_n}&x_2^{\lambda_n}&\cdots&x_n^{\lambda_n}
\end{pmatrix},
\end{displaymath}
where $\lambda_1\geq \lambda_2\geq \cdots\geq \lambda_n\geq 1$.
Hence, $\lambda=(\lambda_1,\lambda_2,\ldots,\lambda_n)$ is a partition.
 We see that $M$ is uniquely determined by this partition

We thus connect $M$ with the Schur polynomial $s_\lambda(x_1,x_2,\ldots,x_n)$. In this way all proved results becomes identities for Schur polynomials. Consider the following partition
\begin{displaymath}\{j_1,j_2,\ldots,j_r\}=\{1,2,\ldots,n+r\}\setminus\{\lambda_n,\lambda_{n-1}+1,\ldots,\lambda_1+n-1\}.
\end{displaymath}
We see that the set $\{j_1,j_2,\ldots,j_r\}$ is uniquely determined by $\lambda$ and
$r$.
Finally, We denote $nr+j_1+\cdots+j_r+\frac{r(r-1)}{2}+\lfloor\frac{n}{2}\rfloor
=\nu(\lambda,r)$.

Hence, the following formula is true:
\begin{displaymath}
s_\lambda(x_1,\ldots,x_n)=(-1)^{\nu(\lambda,r)}\cdot\det Q(\lambda,r).
\end{displaymath}

The expression on the left-hand side of this equation is the Schur polynomial
$s_\lambda(x_1,x_2,\ldots,x_n)$.

\begin{remark}Hence, each Schur polyinomial is, up to the sign, a sub-determinant of order $r$ of $P$.
\end{remark}
\begin{remark}
We obtained formula in which the Schur polynomials are obtained in terms of
elementary symmetric polynomials.
Such a formula is  known as the second Jacobi-Trudi identity.
\end{remark}

\end{document}